
\magnification 1200
\font\abs=cmr8
\font\refit=cmti8
\font\refbf=cmbx8

\font\refma=cmmi9
\font\ccc=cmcsc10

\def\cor #1#2 {{\ccc (#1.#2) Corollary.}\phantom{X}}
\def\d{{\Delta}}

\def\ddv{{\fraz{\partial}{\partial v}\,}}

\def\dim{{\sl Proof.}\phantom{X}}
\def\disp#1{{\displaystyle{#1}}}
\def\dfn #1#2 {{\ccc (#1.#2) Definition.}\phantom{X}}
\def\esp#1{e^{\displaystyle#1}}
\def\f{{\cal F}}
\def\ff#1{{\cal F}(#1)}
\def\ffq#1{{\cal F}_q(#1)}
\def\fidi{\hskip5pt \vrule height4pt width4pt depth0pt}
\def\fmu{\phi_{m,u}}
\def\fmus{\phi^*_{m,u}}
\def\fq{{\cal F}_q}
\def\fraz#1#2{{\strut\displaystyle #1\over\displaystyle #2}}
\def\g{{G(1)}}
\def\h{{H}}
\def\hc#1{{{\cal H}_{#1}}}
\def\ho{{H}_0}
\def\hirr{{H}_{m,u}^{\,{\rm irr}}}
\def\hoirr{{H}_{0}^{\,{\rm irr}}}
\def\hh#1{{H_{#1}}}
\def\id{{\rm id}}
\def\ii#1{\item{\phantom{1}#1 .\phantom{x}}}
\def\iii#1{\item{\phantom{1}#1}}
\def\ind#1#2#3{{\rm Ind}(#1\!\!\uparrow \!#2,\,#3)}

\def\lem #1#2 {{\ccc (#1.#2) Lemma.}\phantom{X}}
\def\rem #1#2 {{\ccc (#1.#2) Remark.}\phantom{X}}
\def\ll#1{{Lie(#1)}}

\def\muh{{\widehat{\mu}}}
\def\o#1#2{{{\cal O}_{#1,#2}}}
\def\ome#1{{\omega_{#1}}}

\def\pro #1#2 {{\ccc (#1.#2) Proposition.}\phantom{X}}
\def\ps#1#2{{\Psi_{#1,#2}}}
\def\rmk #1#2 {{\ccc (#1.#2) Remark.}\phantom{X}}
\def\rmks #1#2 {{\ccc (#1.#2) Remarks.}\phantom{X}}
\def\rr#1{{\rho_{#1}}}
\def\rsu#1#2{{\rho^{\uparrow}_{#1,#2}}}

\def\tens{{\,\otimes\,}}
\def\th{{\widehat{t}}}
\def\thm #1#2 {{\ccc (#1.#2) Theorem.}\phantom{X}}

\def\uno{{\bf 1}}
\def\uq{{\cal U}_q}

\def\uuq#1{{\cal U}_q(#1)}
\def\xh{{\widehat{x}}}

\hsize= 15 truecm
\vsize= 22 truecm
\hoffset= 1. truecm
\voffset= 0.3 truecm

\baselineskip= 12 pt
\footline={\hss\tenrm\folio\hss} \pageno=1

\vglue 3truecm
\centerline{\bf INDUCED REPRESENTATIONS OF THE ONE DIMENSIONAL}
\smallskip
\centerline{\bf QUANTUM GALILEI GROUP.}
\bigskip
\bigskip
\centerline{{\it
F.Bonechi, R. Giachetti, E. Sorace, M. Tarlini,}}
\medskip
\baselineskip= 10 pt

\centerline{{\abs  INFN Sezione di Firenze and}}
\centerline{{\abs Dipartimento di Fisica, Universit\`a di Firenze, Italy. }}

\bigskip
\bigskip
{\refbf Abstract.} {\abs We study the representations of the quantum
Galilei group by a suitable generalization of the Kirillov method on
spaces of non commutative functions. On these spaces we determine a
quasi-invariant measure with respect to the action of the quantum
group by which we discuss unitary and irreducible representations. The
latter are equivalent to representations on ${\refma{\ell^2}}$,
{\refit i.e.} on the space of square summable functions on a one
dimensional lattice.}

\bigskip
\bigskip
\baselineskip=14pt
\noindent {\bf 1.}
In some past works we used the quantum Galilei group $\uuq\g$ to
describe the kinematical symmetry of a spin chain on a one dimensional
lattice [1]. The definition of its algebra follows from the relations
of the subgroup of the Heisen\-berg group with discrete space
translations and from the addition of a generator for time
translations, {\it i.e.} the energy. It turns out that the eigenvalue
equation for the Casimir of this algebra is the Schr\"odinger equation
for the free particle on a one dimensional lattice. The physical
requirement for the energy to be additive produces as a consequence a
non cocommutative coalgebra and, by duality, a non commutative algebra
for the representative functions. The deformation parameter of the
quantum group thereby obtained has the interpreta\-tion of a lattice
spacing, whose vanishing reproduces the Galilei Lie group $\g$ as
symmetry of the continu\-um. The main purpose of this letter is to
make evident the connection between discrete\-ness and non
commutativity in a physically meaningful example, by study\-ing  the
representa\-tions of the quantum Galilei group in their natural
framework of non commutative spaces. This method has already proved to
be quite suitable for the study of harmonic analysis and special
functions on homogeneous spaces [2]. We shall extend in the quantum
group framework the technique of induced representations that gives a
complete description of unitary representa\-tions for nilpotent Lie
groups according to the Kirillov theory. For a meaningful physical
interpretation of the algebra, an involution with non standard
properties is re\-quired: this, in turn, implies a multiplicative
involution on the space of quantized functions where we are going to
represent the quantum group. We shall first construct the actions on
the space of formal series and later on we determine appropriate
subspaces where we define a quasi-invariant measure that permits to
discuss unitary and irreducible representations. We finally prove that
the space obtained in this non commutative context is isometric with
the Hilbert space of square summable functions on the integers ({\it
i.e.} with $\ell^2\,$) and we also remark that we can extend the
Kirillov theorem on the equivalence of representations.
\bigskip
{\bf 2.} The quantum Galilei  group can be obtained by deforming the
algebra $\f$ of the functions on $\g$ by means of a non trivial
1-cocycle $\eta$ with values in $\bigwedge^2\ll\g$ that defines a
Lie-Poisson structure. Its description is as follows. Let
$\{\mu,x,t,v\}$ be a set of canonical coordinates of the second kind
on $\g$; the fundamental fields on the group manifold corresponding to
the generators of $\ll\g$ are
$$
\eqalign{
{}&{\cal X}_M = i\,\fraz \partial{\partial\mu}\,,
\phantom{{\cal X}_T = i\,\fraz \partial{\partial t}\,,}
{\cal X}_P = i\,\fraz \partial{\partial x}\,,\cr
{}&{\cal X}_T = i\,\fraz \partial{\partial t}\,,
\phantom{{\cal X}_M = i\,\fraz \partial{\partial\mu}\,,}
{\cal X}_B =
ix\,\fraz \partial{\partial\mu}+it\,\fraz \partial{\partial x}
+i\,\fraz \partial{\partial v}\,.\cr}
$$
We have:
\medskip
\pro 21 {\it Let $g\equiv(\mu,x,t,v)\in\g$.
The map $\eta : \g \rightarrow \bigwedge^2\ll\g$ given by
$$
\eta(g)=(2a\mu-2axv+atv^2)\,{\cal X}_M\wedge {\cal X}_P
-av^2\,{\cal X}_M\wedge {\cal X}_B - 2av\,
{\cal X}_P\wedge {\cal X}_B
$$
satisfies $\eta(gg')=\eta(g)+{\rm Ad}_{g}\eta(g')$ and therefore is a
(non trivial) $1$-cocycle that defines the Lie-Poisson brackets}
$$
\{\mu,x\}=-2a\mu\,,~~~~~\{\mu,v\}=av^2\,,~~~~~\{x,v\}=2av\,,~~~~~
\{t,\cdot\}=0\,.~~~\fidi
$$
\medskip
Denote by $\fq:=\fq(\g)$ the Hopf algebra generated by $\{\mu,x,t,v\}$,
whose relations are
$$
[\mu,x]=-2a\mu\,,~~~~~~~[\mu,v]=av^2\,,~~~~~~~[x,v]=2av\,,
$$
$t$ being a central element. Then $\fq$ is a
deformation of $\f$ in the direction of the Poisson brackets (2.1) and we
obtain a Hopf algebra if we define the same comultiplication and antipode
as in the classical case:
$$
\eqalign{
{}&\d\mu=\mu\tens\uno+\uno\tens\mu+v\tens x +(1/2)v^2\tens t\,,
\phantom{xxxx}
\d t=t\tens\uno+\uno\tens t\,,\cr
{}&\d x=x\tens\uno+\uno\tens x+v\tens t\,,
\phantom{\mu\mu\mu+(1/2)v^2\tens xx}
\d v =v\tens\uno+\uno\tens v\,.\cr}
$$
$$
S(\mu)=-\mu+vx-(1/2)v^2t\,,~~~~~~~~ S(x)=-x+tv\,,~~~~~~~S(t)=-t\,,
~~~~~~~S(v)=-v\,.
$$

The quantum enveloping algebra $\uq:=\uq(\g)$ is expressed in terms of
the generators $\{M,K,K^{-1},T,B\}$ with relations

$$
[K,T]=0\,,~~~~~~K\,B\,K^{-1}=B+aM\,,~~~~~~[B,T]=\fraz{K-K^{-1}}{2a}\,,
$$
where $M$ is central and $KK^{-1}=K^{-1}K=1$. The Casimir reads
$$
C= MT+\fraz 1{2\,a^2}\,[K+K^{-1}-2]\ . \eqno(2.2)
$$
The coalgebra structure is given by
$$
\eqalign{
{}&\d M=M\tens K+K^{-1}\tens M\,,
\phantom{T=T\uno +\uno T}
\d K=K\tens K\,,\cr
{}&\d T=T\tens\uno +\uno\tens T\,,
\phantom{M=M K+K^{-1}M}
\d B=B\tens K+K^{-1}\tens B\,\cr}
$$
and
$$
S(M)=-M\,,~~~~~~~S(K)=K^{-1}\,,~~~~~~~S(T)=-T\,,~~~~~~~S(B)=-B-aM\,.
$$

There is a nondegenerate duality pairing between $\fq$ and $\uq$ given by
$$
\langle \mu^\alpha x^\beta t^\gamma v^\delta\, ,
I^{\alpha'} K^{\ell} T^{\gamma'} N^{\delta'}  \rangle =
i^{\alpha+\beta+\gamma+\delta}\ \alpha!\,\gamma!\,\delta!\,
(a\,\ell)^\beta\,\,\delta_{\alpha,\alpha'}\,
\delta_{\gamma,\gamma'}\,\delta_{\delta,\delta'}
$$
where $I=K^{-1} M,\, N=K B$ and $\ell\in {\bf Z}$ while the other
indices are in ${\bf N}$.

To define the classical limit we write $K=\esp{iaP}$ and the Casimir
(2.2) for $a\rightarrow 0$ reproduces the classical quadratic form
$C_{{\rm cl}}=MT-P^2/2M$.

In order to have a meaningful physical interpretation of $\uq$ as a
deformation of the classical symmetry, we shall define the following
involution $M^*=M,\,K^*=K^{-1},\,T^*=T,\,B^*=B$. The compatibility
with the Hopf structure results then in
$$
(XY)^*=Y^*\,X^*\,,~~~~~~~~~~\d (X^*)= (\ast\tens\ast)\circ\tau\circ\d X\,,
$$
where $\tau(X\tens Y)=Y\tens X$. In the classification of
$\ast$-structures given in [3], this is referred to as to the case
(1,1). For this type of involution we have $\ast\circ S=S\circ\ast$
and we see that $(\ast\circ S)^2$ is different from identity.
Therefore the use of the duality relations and of the operator
$\ast\circ S$ for determining an involution on $\fq$ requires some
modifications. A direct calculation shows:
\medskip
\lem 23 {\it Let $Q:\uq\rightarrow\uq$ be the morphism defined on the
generators by
$$
Q(M)=M\,,~~~~~Q(K)=K\,,~~~~~Q(T)=T\,,~~~~~Q(B)=B+aM\,.
$$
Then $Q$ is invertible and $Q^2=(\ast\circ S)^2$. Therefore
$\sigma=Q^{-1}\circ\ast\circ S$ is an antilinear multiplicative map
of $\uq$ with $\sigma^2=\id$.
Moreover $\d\circ\sigma=(\sigma\tens\sigma)\circ\d$.}\fidi
\medskip
The involution on $\fq$, again denoted by $\ast$, can now be defined
as follows.
\medskip
\pro 24 {\it Let $f^*$ be defined by the relation
$$
\langle\,f^*,X\,\rangle = \overline
{\langle\,f,\sigma(X)\,\rangle}\,.
$$
Then the map $f\mapsto f^*$ is an involution on $\fq$.

The generators $\{\mu,x,t,v\}$ are real under $\ast$ and we have:
\smallskip
\iii{$\phantom{i}(i)$} $(f_1f_2)^*=f_1^*f_2^*\,;$
\iii{$(ii)$} $\d(f^*)=(\ast\tens\ast)\,\d f\,.$
\smallskip
According to {\rm [3]} the involution is therefore of the
$(0,0)$-type.}
\smallskip
\dim Since $\sigma$ and complex conjugation are
commuting involutions, it is immediate to see that also
$f\mapsto f^*$
is an involution.

For the point $(i)$ we have:
$$
\eqalign{
{}\langle (f_1f_2)^*,X\rangle&=\overline{\langle f_1\tens f_2,
\d\circ\sigma(X)\rangle} \cr
{}&=\sum\limits_{(X)}\overline{\langle f_1\tens f_2,
 \sigma(X_{(1)})\tens\sigma(X_{(2)})\rangle}=
\langle f_1^*f_2^*,X\rangle\cr
}
$$
For the second item:
$$
\eqalign{
{}\langle \d f^*,X\tens Y\rangle&=\overline{\langle f,
\sigma(XY)\rangle} \cr
{}&=\sum\limits_{(f)}\,\overline{\langle f_{(1)},
\sigma(X)\rangle}\,\,\overline{\langle f_{(2)},
\sigma(Y)\rangle}=
\langle (\ast\tens\ast)\d f,X\tens Y\rangle\cr
}
$$
Finally, as $\sigma(I^\alpha K^\ell T^\gamma N^\delta)
=(-)^{\alpha+\gamma+\delta}\,I^\alpha K^\ell T^\gamma N^\delta$,
the generators of $\fq$ are real.
\fidi

\bigskip
{\bf 3.} In this section we present some algebraic aspects of induced
representations that will be useful for applications to the quantum
Galilei group. We take ad\-van\-tage from the fact that the undeformed
group $\g$ is nilpotent, so that its representations are completely
described by Kirillov theory [4]. In particular, due to the existence
of maximal polarizing subalgebras, any irreducible representation is
always induced from a character $\ome J$ of a subgroup $J$. The
representative space $\hh \omega$ is the space of complex valued functions
that are $\ome J$-covariant along the $J$-cosets,
$$
f(jg)=\ome J(j) f(g)\,,~~~~~~~~~~~~j\in J\,,\, g\in G\,,
\eqno(3.1)
$$
and square integrable with respect to an invariant measure $\nu$ on
the homoge\-ne\-ous space $J\!\setminus\!G$. Therefore $\hh \omega$ is
isometric to $L^2(J\!\setminus\!G,\nu)$ and $\rr
\omega=\ind{J}{\g}{\ome J}$ reduces to the restriction to $\hh \omega$
of the regular represen\-ta\-tion
$$
\Bigl(\rr \omega(g_2)f\Bigr)(g_1)=f(g_1g_2)\,,~~~~~~~~~~~~g_1,g_2\in
G\,. \eqno(3.2)
$$

We now reformulate this scheme in a Hopf algebra framework. Let us
examine first the covariance condition (3.1). Observe that $\ff J$ is
a Hopf algebra and $\pi_J: \ff{G}\rightarrow \ff{J}$, defined by
$\pi_J(f)=f|_J$, is a Hopf algebra homomorphism. The character $\ome
J$ determines a corepresentation $1\rightarrow \ome J \,:\,{\bf
C}\rightarrow {\bf C}\tens\ff{J}\simeq\ff{J}$. Consider the right
action of $j\in J$ on $f\in\f$ given by $ (f\cdot j)(g)=f(jg)$ and let
$\Lambda:\ff{G}\rightarrow\ff{J}\tens\ff{G}$ be the corresponding left
coaction. The actual form of $\Lambda$ is given by the standard
dualization procedure, namely $\Lambda=(\pi_{J}\tens\id)\circ\d$. The
space $\hh \omega$ of the induced representation will be defined as
the subspace of functions in $\ff{G}$ that satisfy the equivariance
condition (3.1) rewritten as
$$
\Lambda(f)=\ome J\tens f\ , ~~~~~~~~~~~~~~f\in \ff{G}\ . \eqno(3.3)
$$
On $\hh \omega$ the comultiplication $\d$ determines a
corepresentation $\Psi_\omega$.

The generalization of the procedure to quantum groups is
straightforward when a quantum subgroup exists, as in this case. We
define the Hopf algebra $\ffq J$ generated by the three primitive
elements $\{\muh,\xh,\th\,\}$ with relations $[\muh,\xh]=-2a\muh$ and
$[\th,\cdot]=0$. We also assume that the involution is of the type
$(0,0)$ and that these three generators are real. The map
$\pi_J:\fq\rightarrow\ffq J$, defined on the generators by
$\pi_J(\mu)=\muh,\,\pi_J(x)=\xh,\,\pi_J(t)=\th,\,\pi_J(v)=0$, is a
surjective $\ast$-Hopf morphism. It is easy to verify that
$\ome{m,u}=\exp[- i(m\muh+u\th\,)]$ defines a one dimensional
corepresentation of $\ffq J$.

\medskip
\pro 34 {\it Let $\hc{m,u} =\Bigl\{\fmu\,f(v)\Bigr\}$, where
$\fmu=\exp[- i(m\mu+ut)]$ and $f$ is a formal series in $v$.
\iii {$\phantom{ii}(i)$} Each element $\fmu\,f(v)\in\hc{m,u}$ is
equivariant according to $(3.3)$;
\iii {$\phantom{i}(ii)$} on  $\hc{m,u}$ we can define
a coaction by
$$
\eqalign{
\ps mu\Bigl(\fmu\, f(v)\Bigr)=&\exp[-im(\mu\tens\uno+\uno\tens\mu+
v\tens x+(1/2)\,v^2\tens t)]\cr
{}&\phantom{XXXXXX}\cdot\,\exp[-iu(t\tens\uno+\uno\tens t)]\,
\,f(v\tens\uno+\uno\tens v)\,.\cr}
$$
The corresponding action of an element $X\in \uq$ is given by
$\rr{m,u}(X)=(\id\tens X)
\ps mu$;
\iii {$(iii)$} letting $\rsu mu(X) = \phi_{m,u}^{-1}\,\,\rr{m,u}(X)\,\,
\fmu$, with $X\in\uq$, we have}
$$
\eqalign{
{}&\rsu mu(M)\,f(v) = m\,f(v)\,,
\phantom{K = \fraz 1{1-iamv}}
\rsu mu(T)\,f(v) = \Bigl[\fraz{mv^2}{2(1-iamv)}+u\Bigr]
\,f(v)\,,\cr
{}&\rsu mu(K)\,f(v) = \fraz 1{1-iamv}\,f(v)\,,
\phantom{M = m}
\rsu mu(B)\,f(v) = i\,(1-iamv)\,\ddv f(v)\,.\cr}
$$

\medskip
For the proof we need the following simple technical results.
\medskip
\lem 35 {\it For complex numbers $r,s$ and natural number $n$,
the following holds:}
\smallskip
\iii {$\phantom{ii}(i)~~$} $\esp{r(\mu-s\,v)}=\esp{r\mu}\,
     (1+ar\,v)^{\disp{-s/a}}
       =(1-ar\,v)^{\disp{s/a}}\,\esp{r\mu} \,;$
\iii  {$\phantom{i}(ii)~~$} $\esp{r(\mu-s\,v^2)}=\esp{r\mu}\,
        \esp{-(srv^2(1+arv)^{-1})}\,;$
\iii {$(iii)~~$} $\mu^{\disp{n-1}}\,(\mu+na\,v)=(\mu+a\,v)^{\disp n}\,.$
\smallskip
\dim $(i)$ Suppose that $L(v)$, as a formal series of $v$, satisfies the
relation $[L(v),\mu] =s\,v$. Then it is easily verified that
$$
\esp{r(\mu-s\,v)}=\esp{-L(v)}\,\esp{r\mu}\,\esp{L(v)}
$$
Therefore, since
$\exp[-L(v)]\,\exp[r\mu]=\exp[r\mu]\,\exp[-L(v)]\,(1+ar\,v)^
{\displaystyle{-s/a}}$,
the final result follows by substituting this expression into the
former one. With a similar method also item $(ii)$ can be proved. The
last statement follows by an easy induction. \fidi

\medskip
{\sl Proof of Proposition $(3.4)$} The proof is a consequence of the
duality relations, the properties of the comultiplication and of the
results of Lemma (3.5). We report here only on some more meaningful
points. For getting the action of the mass, a first computation shows
that
$$
(\id\tens M)\,(\d\mu)^n=in\mu^{n-1}\,.
$$
It is then straightforward to find the result.

Concerning $K$ a useful relation turns out to be
$$
(\id\tens K)\,(\d\mu)^{\disp n}=(\mu-av)^{\disp n}\,,
$$
as a consequence of which we find
$$
\rr{m,u}(K)\Bigl(\fmu\,f(v)\Bigr)=\esp{-im(\mu-a\,v)}\,\esp{-iut}
\,f(v)\,.
$$
Using $(3.5i)$ to factorize $\fmu$ we find $\rsu mu(K)$ as in
the proposition.

In order to obtain the action of the energy, we point out that
$$
(\id\tens \esp{irT})\,(\d\mu)^{\disp n}=(\mu-rv^2/2)^{\disp n}
$$
and then, by differentiating with respect to $r$ and using $(3.5ii)$,
we get the result. Finally, the boost does not present any particular
problem.
\fidi
\bigskip
{\bf 4.} The framework of formal series is not suitable to discuss
unitary represen\-ta\-tions. In order to determine an inner product, a
quasi-invariant measure has to be defined on the representation space.
We must therefore find appropriate subspaces of formal series that
carry such a measure. Although we will not be able to define a
coaction, we show that an action is well defined on these
restrictions: this is sufficient to define representations that will
prove to be unitary. Irreducibility is finally discussed.
\medskip
\lem 41 {\it $($i$\,)$ The following holds
 $$(1+iamv)\,\fmu\,v\,\fmu^{-1}=v\,;$$

$($ii$\,)$ let $v_k=(\fmu)^{-k}\,v\,(\fmu)^k\,.$ Then
$$
iam(k-\ell)v_k\,v_\ell=v_k-v_\ell\ .~~~~~\fidi
$$}

Consider the linear space freely generated over the set of finite
words in $\fmu\,,\,\fmu^{-1}\,,$ and $v$ with relations (4.1$i$). Let
$\ho$ be the subspace generated over the words whose total degree in
$\fmu$ is equal to that in $\fmu^{-1}$.
\medskip
\pro 42 {\it If $\ho^{(k)}=v_k\,{\bf C}[v_k]$, i.e. the space of
polynomials in
$v_k$ with vanishing constant term, then}
$$\ho={\bf C}\oplus \Bigl(\bigoplus\limits_{k=-\infty}^{\infty}\ho^{(k)}
\Bigr)\,.$$
\smallskip
\dim It is a direct consequence of the definitions and of relation 
$(4.1ii)$. \fidi
\medskip

\medskip
Let $\h_{m,u}^{(k)}=\fmu\,\ho^{(k)}$ and $\h_{m,u}=\fmu\,\ho$. Since
the coproduct of $\fmu$ is not polynomial in $\fmu\,,\fmu^{-1}\,,$ and
$v$ no coaction of the quantum Galilei group can be defined on $\ho$
and $\h_{m,u}$. However the action obtained as in $(3.4ii)$ admits a
restriction to $\h_{m,u}$ and  therefore determines a
represen\-ta\-tion of $\uq$.
\medskip
\pro 43 {\it $\uq$ is represented on $\h_{m,u}$ as follows:
$$
\eqalign{
{}&\rr{m,u}(M)\,\fmu\,v_k^n=m\fmu\,v_k^n\ ,\cr
{}&\rr{m,u}(K)\,\fmu\,v_k^n=\fmu\,v_k^n(1+iam v_1)\ ,\cr
{}&\rr{m,u}(K^{-1})\,\fmu\,v_k^n=\fmu\,v_k^n (1-iam v_0)\ ,\cr
{}&\rr{m,u}(T)\,\fmu\,v_k^n=\fmu\,v_k^n\,
                    (\,\fraz i{2a}(v_0-v_1)+u)\ , \cr
{}&\rr{m,u}(B)\,\fmu\,v_k^n=i\,n\,\fmu v_k^{n-1}\,(1-imv_0)
(1+iam\,k\,v_k)^2\ . \cr
}
$$
\noindent
The representation is reducible but not completely reducible. The only
irreducible component is the restriction to
$$\hirr={\bf C}\oplus\h_{m,u}^{(0)}\oplus
\h_{m,u}^{(1)}\,.$$}
\smallskip
\dim The proof is obtained by a direct calculation from $(3.4iii)$ and
$(4.1ii)$.
\fidi
\medskip
\cor 44 {\it
Let $\chi_m=(1+iamv_{1})\in{\bf C}\oplus
\ho^{(1)}$. This element is invertible in $\ho$ with $\chi_m^{-1}=(1-iamv_0)
\in{\bf C}\oplus \ho^{(0)}$. Then $\{\fmu\,\chi_m^\ell\}_{\ell\in {\bf
Z}}$ is a basis for $\hirr$ and the irreducible representations given
in $(4.3)$ read
$$
\eqalign{
{}&\rr{mu}(K^{\pm 1})\,\fmu\, \chi_m^\ell= \fmu\, \chi_m^{\ell\pm 1}\,,
\quad\quad
\rr{mu}(B)\,\fmu\, \chi_m^\ell = m\,\ell\,\fmu\, \chi_m^\ell\,,\cr
{}&\rr{mu}(T)\,\fmu\, \chi_m^\ell = \fmu\, \chi_m^\ell\, (\,\fraz
1{2a^2m}\,(2-\chi_m-\chi_m^{-1})+u)\,, }
$$
while $\rr{mu}(M)$ is the multiplication by $m$ }.\fidi
\medskip
Let $\hoirr={\bf C}\oplus \ho^{(0)}
\oplus \ho^{(1)}$. Define a linear
functional $\nu_a: \hoirr\rightarrow {\bf C}$ by
$$
\nu_a(1)= 1\ , \quad\quad
\nu_a(v_0^n)=\fraz 1{(iam)^n}\ , \quad\quad \nu_a(v_{1}^n)=\fraz
1{(-iam)^n}\ , \eqno{(4.5)}
$$
so that $\nu_a(\chi_m^\ell)=\delta_{\ell,0}$.
\medskip
\pro 46 {\it The following holds:

\smallskip
\iii {$\phantom{ii}(i)$} if $a,b\in \hirr$ then
$a^*b\in \hoirr$ and $\langle a,b\rangle=\nu_a(a^*b)$ defines a scalar
product on $\hirr\,$;

\smallskip
\iii {$\phantom{ii}(ii)$} the Hilbert space obtained by completing
$\hirr$ with respect to the scalar prod\-uct defined in $(i)$ is
isometric to $L^2(\,[0,2\pi/a)\,,\,a\,dp/(2\pi))\,$ and, by Fourier
transform, to $\ell^2$;

\smallskip
\iii {$\phantom{ii}(iii)$} the representation in $(4.4)$ is unitary.}
\smallskip
\dim  The first statement of ($i$) is a direct
consequence of the relation $\fmus \chi_m^*\fmu=\chi_m^{-1}$.
Moreover
$$
\langle\fmu\,\chi_m^n,\fmu\,\chi_m^\ell\rangle =
\nu_a(\fmus(\chi_m^n)^*
\fmu\,\chi_m^\ell)=\nu_a(\chi_m^{\ell-n})=\delta_{\ell,n}\ .
$$
Thus
$$
\langle\fmu\,\chi_m^n,\fmu\,\chi_m^\ell\rangle=\fraz a{2\pi}\int_0^{2\pi/a}
\esp{iap(\ell-n)}dp\,.
$$
We therefore obtain a natural isometry of $\hoirr$ into
$L^2(\,[0,2\pi/a)\,,\,a\,dp/(2\pi))$  map\-ping $\chi_m^\ell$ in plane
waves.

Item ({\it iii}) is proved by verifying that
$\rr{m,u}(X^*)=\rr{m,u}(X)^t$ on the generators.\fidi
\medskip
\cor 47 {\it The measure $\nu_a$ is quasi-invariant with respect to
the restriction to $\hoirr$ of the regular representation of $\,\uq$. 
Indeed for each $f\in \hoirr$ and $X\in\uq$ we have
$$
\nu_a((\id \tens X)\d f) = \nu_a(f\,\xi(X))\,
$$
where}
$$
\xi\mapsto\xi(X)=\sum_{n=0}^\infty (iam)^n\langle v^n,X\rangle\,\chi^n_m\,
:\,\uq\rightarrow\hoirr\,.~~~\fidi
$$
\medskip
\rmks 48
We shall conclude by pointing out some peculiar differences between
the classical and the quantum situation and by describing the
extension of the Kirillov theorem to the quantum Galilei group.
\smallskip
($i$) The measure $\nu_a$ is finite on the coordinate $v$, at
difference with  the classical case where it diverges, as can be
checked by taking the limit $a\rightarrow 0$. Noncommuta\-tivity means
here regularization.

($ii$) Contrary to the classical situation, $\nu_a$ is only
quasi-invariant and becomes invariant only in the continuum limit,
[5]. Nevertheless these induced represen\-ta\-tions are unitary. This
feature is closely connected with the non standard properties of the
involution.

($iii$) On $\hoirr$ we have a scalar product defined by $\langle
a,b\rangle
= \nu_a(a^\dagger b)$ where $a^\dagger=\fmus a^*\fmu$. The Hilbert
space obtained is isometric with the space given in (4.6). If we
formally write $\chi_m=\sum_{k=0}^\infty(1-\chi_m^{-1})^k$, we note
that this series doesn't converge in norm.

($iv$) The classical limit of the representation (4.4)  corresponds to
the choice of the point $mM^*+uT^*$ on the coadjoint orbit $\o mu
=\Bigl\{mM^*+\alpha P^* + (u+\alpha^2/2m)\, T^* + \beta B^*
\,|\, \alpha,\beta \in {\bf R} \Bigr\}\,,$ characterized
by the mass $m$ and by the ``internal energy'' $u$. A generic point of
the orbit $\o mu$ gives a character
$\exp[-i(m\muh+\alpha\xh+(u+\alpha^2/2m)\th\,)]$ that  determines a
representation $\rho$ on the space of equivariant functions
$H=\{\phi\,f(v)\}$, where $\phi=\exp[-i(m\mu+\alpha
x+(u+\alpha^2/2m)t)]$. This classical procedure extends to the quantum
case. For each triple of numbers ${\bf c}=(c_1,c_2,c_3)$ we have that
$\omega_{\bf c}=\exp[-i c_1\muh]\cdot
\exp[-i c_2\th\,]\cdot \exp[-i c_3 \xh]$ defines a one dimensional
corepresentation of $\ffq J$. The corresponding induced action
$\rho^\uparrow = \phi_{\bf c}^{-1}\,\rho\,
\phi_{\bf c}$ is
$$
\eqalign{
\rho^\uparrow (M)\,f(v) &= c_1\exp(ia\,c_3)\,f(v)\,,\cr
\rho^\uparrow (K)\,f(v) &=
(\exp(-ia\,c_3)-iac_1\exp(ia\,c_3)\,v)^{-1}\,f(v)\,,\cr
\rho^\uparrow (B)\,f(v) &= i\,(\exp(-ia\,c_3)-ia\,c_1\exp(ia\,c_3)\,v)\,
\ddv f(v)\,,\cr
\rho^\uparrow (T)\,f(v) &=
\Bigl(\fraz{c_1\exp(4ia\,c_3)\,v^2}{2(1-ia\,c_1\exp(2ia\,c_3)\,v)}-
\fraz{1-\exp(2ia\,c_3)}{2ia}\,v + c_2\Bigr)\,f(v)\,,\cr}
$$
with $\phi_{\bf c}= \exp[-i c_1\mu]\cdot
\exp[-i c_2\,t]\cdot \exp[-i c_3\,x]$.

The calculation of the Casimir (2.2) gives
$$
C=c_1\,c_2\,\esp{iac_3}+\fraz1{a^2}(\cos(ac_3)-1)\,,
$$
so that the equivalence with the representation in $(3.4iii)$ is
proved posing
$$
c_1=m\,\esp{-iac_3}\quad\quad {\rm and} \quad\quad
c_2=-\fraz1{ma^2}(\cos(ac_3)-1)+ u\,.
$$

The construction of unitary representations requires that $c_3$ is
real, so that $\omega_{\bf c}^*\omega_{\bf c}=1$ and the explicit
equivalence is easily obtained by defining the plane waves as $\
\chi_{m,c_3}=\esp{iac_3}+iam\,\phi^*_{\bf c}\,v\,\phi_{\bf c}\ $ and
$\ \chi^{-1}_{m,c_3}=\esp{-iac_3}-iam\,v\ $.

\vfill\break
\baselineskip= 12 pt
\bigskip
\bigskip

\centerline{{\bf References.}}

\bigskip
\baselineskip= 10 pt
{\abs
\ii {1}
Bonechi F., Celeghini E., Giachetti R., Sorace E. and Tarlini M.,
      {\refit Phys. Rev. B} {\refbf 32}, 5727 (1992) and
      {\refit J. Phys. A} {\refbf 25}, L939 (1992).
\smallskip
\ii {2}
        Bonechi F., Ciccoli N., Giachetti R., Sorace E. and Tarlini M.,
      {\refit Commun. Math. Phys.} {\refbf 175}, 161 (1996).
\smallskip
\ii {3}
       Scheunert M., {\refit J. Math. Phys.} {\refbf 34}, 320 (1991).
\smallskip
\ii {4}
        Corwin L. and  Greanleaf F.P. ``{\refit Representations of
     Nilpotent Lie Groups and their Applica\-tions}'', (Cambridge
     University Press, Cambridge, 1990).
\smallskip
\ii {5}
       Kirillov A.A., ``{\refit Elements of the Theory of
      Representations}'', (Springer Verlag, Berlin, 1990).
\smallskip

}
\bye